\documentclass{amsart}

\usepackage{amsmath,amssymb,amsthm}

\theoremstyle{definition}
\newtheorem{definition}{Definition}

\theoremstyle{plain}

\newtheorem{theorem}[definition]{Theorem}
\newtheorem{corollary}[definition]{Corollary}

\begin{document}

\title[How to compute the Wedderburn decomposition]
{How to compute the Wedderburn decomposition of a finite-dimensional
associative algebra}

\author{Murray R. Bremner}

\address{Department of Mathematics and Statistics, University of Saskatchewan,
Canada}

\email{bremner@math.usask.ca}



\begin{abstract}
This is a survey paper on algorithms that have been developed during the last
25 years for the explicit computation of the structure of an associative
algebra of finite dimension over either a finite field or an algebraic number
field. This constructive approach was initiated in 1985 by Friedl and R\'onyai
and has since been developed by Cohen, de Graaf, Eberly, Giesbrecht, Ivanyos,
K\"uronya and Wales. I illustrate these algorithms with the case $n = 2$ of the
rational semigroup algebra of the partial transformation semigroup $PT_n$ on
$n$ elements; this generalizes the full transformation semigroup and the
symmetric inverse semigroup, and these generalize the symmetric group $S_n$.
\end{abstract}

\maketitle


\section*{Introduction}

Part 1 of this survey begins by recalling the classical structure theory of
finite-dimensional associative algebras over a field; the most important
results are Dickson's Theorem characterizing the radical in characteristic 0,
the Wedderburn-Artin Theorem on the structure of semisimple algebras, and the
Wedderburn-Malcev Theorem on lifting the semisimple quotient to a subalgebra.
It continues by quoting observations from Friedl and R\'onyai \cite{FriedlR} to
motivate a constructive computational approach to the theory. This explicit
approach requires a presentation of the algebra by a basis and structure
constants, and algorithms for calculating the following: a basis for the
radical of the algebra; structure constants for the semisimple quotient; a
basis for the center of the semisimple quotient; a new basis for the center
consisting of orthogonal idempotents; the identity matrices in the simple
ideals of the quotient; an isomorphism of each simple ideal with a full matrix
algebra; explicit matrices for the irreducible representations; and a
subalgebra isomorphic to the semisimple quotient. This survey emphasizes
characteristic 0: in this case, all calculations can be reduced to computing
the row canonical form of a matrix.

Part 2 begins by introducing some classical semigroups of Boolean matrices
which are natural generalizations of the symmetric group. The main example is
the semigroup of partial transformations on $n$ elements. It continues by
presenting explicit calculations for $n = 2$ to illustrate the theory and
algorithms of Part 1.


\section{Theory and algorithms}

\subsection{Structure theory of associative algebras}

We consider only associative algebras $A$ of finite dimension over a field $F$.
We usually assume that $F$ is a finite extension of either the field
$\mathbb{Q}$ of rational numbers or the field $\mathbb{F}_p$ with $p$ elements
($p$ prime); that is, an algebraic number field or a finite field. To keep the
exposition as simple as possible, we often assume that $F = \mathbb{Q}$. For
the classical structure theory of finite-dimensional associative algebras, our
main reference is Drozd and Kirichenko \cite{DrozdK}. For an account of the
historical development, see Parshall \cite{Parshall}.

\begin{definition}
\cite[\S 2.2]{DrozdK} A left $A$-module $M$ is \textbf{semisimple} if it is
isomorphic to a direct sum of simple modules. An algebra $A$ is
\textbf{semisimple} if its left regular module is semisimple. A left ideal $I$
of $A$ is \textbf{nilpotent} if $I^m = \{0\}$ for some $m \ge 1$. An element $x
\in A$ is \textbf{strongly nilpotent} if the principal left ideal $Ax$ is
nilpotent.
\end{definition}

\begin{theorem}
\cite[Corollaries 2.2.5, 2.2.6]{DrozdK} The following conditions are
equivalent: ($i$) $A$ is semisimple; ($ii$) $A$ contains no nonzero nilpotent
left ideals; ($iii$) $A$ contains no nonzero strongly nilpotent elements.
\end{theorem}

\begin{theorem}
\cite[Theorem 2.4.3, Corollary 2.4.5]{DrozdK} \emph{(Wedderburn-Artin Theorem)}
Every semisimple algebra $Q$ has a unique decomposition $Q = Q_1 \oplus \cdots
\oplus Q_c$ into the direct sum of simple ideals where $Q_i Q_j = \{0\}$ for $i
\ne j$. Every simple algebra is isomorphic to a full matrix algebra $M_n(D)$
for some division algebra $D$ over $F$.
\end{theorem}

\begin{definition}
\cite[\S 3.1]{DrozdK} The \textbf{radical} $R(M)$ of a left $A$-module $M$
consists of all $y \in M$ such that $f(y) = 0$ for every homomorphism $f$ from
$M$ to a simple left $A$-module. The \textbf{radical} $R(A)$ of the algebra is
the radical of the left regular module.
\end{definition}

\begin{theorem}
\cite[Theorems 3.1.6, 3.1.10]{DrozdK} The radical $R(A)$ is the set of all
strongly nilpotent elements; it is a two-sided ideal and $Q = A / R(A)$ is
semisimple.
\end{theorem}

\begin{definition}
\cite[\S 6.1]{DrozdK} An algebra $A$ over a field $F$ is \textbf{separable} if
the scalar extension $A \otimes_F K$ is semisimple for every field extension
$K$ of $F$.
\end{definition}

\begin{theorem}
\cite[Corollary 6.1.4]{DrozdK} Every separable algebra is semisimple; the
converse holds if $F$ is a perfect field (in particular, if $\mathrm{char}\,F =
0$ or $F$ is finite).
\end{theorem}

\begin{definition}
\cite[\S 6.2]{DrozdK} Let $\pi\colon A \to Q = A / R(A)$ be the canonical
surjection. A \textbf{lifting} of $Q$ to $A$ is a homomorphism $\epsilon\colon
Q \to A$ such that $\pi\epsilon$ is the identity on $Q$.  It is clear that
$\epsilon$ is injective, that $\epsilon(Q)$ is a subalgebra of $A$ isomorphic
to $Q$, and that $A = \epsilon(Q) \oplus R(A)$ as vector spaces. Two liftings
$\epsilon$ and $\eta$ are \textbf{conjugate} if there is an invertible element
$a \in A$ such that $\eta(x) = a^{-1} \epsilon(x) a$ for all $x \in Q$, and
\textbf{unipotently conjugate} if $a = 1 + \zeta$ for some $\zeta \in R(A)$.
\end{definition}

\begin{theorem} \label{wedderburnmalcev}
\cite[Theorem 6.2.1]{DrozdK} \emph{(Wedderburn-Malcev Theorem)} If $Q = A/R(A)$
is separable then a lifting exists and any two liftings are unipotently
conjugate.
\end{theorem}

\subsection{A constructive approach to the classical theory}

As motivation for a computational approach, we quote the following passages
(with slight changes) from Friedl and R\'onyai \cite[\S\S 1.1, 1.2,
1.4]{FriedlR}: ``The textbook proofs of these results are not constructive.
They mostly start by picking `any minimal [left] ideal'. But the minimal [left]
ideals may not cover more than a tiny fragment of the algebra and might be
quite difficult to find. {\dots} Finding the radical and the simple factors of
the [semisimple] quotient are as essential to computational algebra as
factoring integers and finding composition factors are to computational number
theory and group theory. {\dots} Such results are likely to have applications
to computational group theory as well since group representations are a major
source of problems on matrix algebras. {\dots} The case of commutative
associative algebras generalizes the problem of factoring polynomials over [a
field] $F$. Indeed, let $f \in F[x]$ and let $f = g_1^{e_1} \cdots g_k^{e_k}$
where the $g_i$ are irreducible over $F$. Consider the commutative associative
algebra $A = F[x] / \langle f \rangle$. The radical of $A$ comes from the
`degeneracy' of $f$, i.e.~the presence of multiple factors: $R(A)$ is generated
(as an ideal of $A$) by $h = g_1 \cdots g_k$.  The quotient $A/R(A)$ is
isomorphic to $F[x]/\langle h \rangle$. This in turn is the direct sum of its
simple components, the fields $F[x]/\langle g_i \rangle$	($i=1,...,k$).	
Finding these components is equivalent to factoring $f$.''

\subsection{Limitations of this survey}

The goal of this brief survey is to present the essential ideas in enough
detail that the algorithms can be translated more-or-less directly into
computer programs.  Therefore, some important issues are ignored, but
references will be given: ($i$) Computational complexity: Most of the
algorithms terminate in a number of steps which is a polynomial function of the
size of the input. ($ii$) Computing the radical in characteristic $p$: This is
much more difficult than in characteristic 0. ($iii$) The possibility that the
minimal polynomials of central elements do not split over the base field: This
seems like a severe restriction, but it is satisfied by many important
examples, such as the group algebra of the symmetric group. ($iv$) The general
case of finding a minimal left ideal in a simple ideal of the semisimple
quotient: This is equivalent to computing an explicit isomorphism of the simple
ideal with a full matrix algebra.

\subsection{Structure constants}

Since the algebra $A$ is finite dimensional over the field $F$, it is
completely determined by a basis $\{ a_1, \dots, a_n \}$ over $F$ and structure
constants $c_{ij}^k \in F$ such that
  \[
  a_i a_j = \sum_{k=1}^n c_{ij}^k a_k
  \quad
  (1 \le i, j, k \le n).
  \]

\subsection{The radical: Dickson's theorem}

The definition of the radical does not depend on the base field, and so we can
regard $A$ as an algebra over $\mathbb{Q}$ or $\mathbb{F}_p$. In characteristic
0, Dickson's Theorem reduces finding a basis for the radical to solving a
linear system. In characteristic $p$, the problem is more difficult; see Friedl
and R\'onyai \cite{FriedlR}, R\'onyai \cite{Ronyai1}, Cohen et al.
\cite{CohenIW}. In this survey we consider only $F = \mathbb{Q}$.

\begin{definition}
For $x \in A$ the \textbf{left multiplication operator} $L_x \in
\mathrm{End}_F(A)$ is $L_x(y) = xy$, and $[L_x]$ is its matrix with respect to
the given basis of $A$.
\end{definition}

We assume that $A$ is unital, adjoining an identity if necessary; then the
representation $x \mapsto [L_x]$ of $A$ is faithful and $A$ is isomorphic to a
subalgebra of $M_n(F)$.

\begin{theorem} \cite[\S 65]{Dickson} \emph{(Dickson's Theorem)}
If $\mathrm{char}\,F = 0$ and $A$ is a subalgebra of $M_n(F)$ then $x$ is in
the radical of $A$	if and only if $\mathrm{trace}(xy) = 0$ for every $y \in
A$.
\end{theorem}

We use this to express the radical as the nullspace of a matrix. Let $x$ be a
linear combination of $\{ a_1, \dots, a_n \}$ such that $\mathrm{trace}(xy) =
0$ for every $y$. By linearity, it suffices to assume $\mathrm{trace}(xa_i) =
0$ for $i = 1, \dots, n$. For $x_j \in F$ we have
  \allowdisplaybreaks
  \begin{align*}
  &
  x = \sum_{j=1}^n x_j a_j \in A,
  \quad
  x a_i
  =
  \sum_{j=1}^n x_j a_j a_i
  =
  \sum_{j=1}^n x_j \sum_{k=1}^n c_{ji}^k a_k
  =
  \sum_{k=1}^n \sum_{j=1}^n c_{ji}^k x_j a_k,
  \\
  &
  x a_i a_\ell
  =
  \sum_{k=1}^n \sum_{j=1}^n c_{ji}^k x_j a_k a_\ell
  =
  \sum_{k=1}^n \sum_{j=1}^n c_{ji}^k x_j \sum_{m=1}^n c_{k\ell}^m a_m
  =
  \sum_{m=1}^n \sum_{j=1}^n \sum_{k=1}^n c_{ji}^k c_{k\ell}^m x_j a_m.
  \end{align*}
Hence the matrix representing left multiplication by $x a_i$ and its trace are
as follows:
  \[
  [L_{xa_i}]_{m\ell}
  =
  \sum_{j=1}^n \sum_{k=1}^n c_{ji}^k c_{k\ell}^m x_j,
  \quad
  \mathrm{trace}([L_{xa_i}])
  =
  \sum_{j=1}^n \Big( \sum_{k=1}^n \sum_{\ell=1}^n c_{ji}^k c_{k\ell}^\ell \Big) x_j.
  \]

\begin{corollary}
The radical of $A$ is the nullspace of the $n \times n$ matrix $\Delta$ such
that
  \[
  \Delta_{ij} = \sum_{k=1}^n \sum_{\ell=1}^n c_{ji}^k c_{k\ell}^\ell.
  \]
\end{corollary}

If $A$ is a semigroup algebra then $a_i a_j = a_{\mu(i,j)}$ and $c_{ij}^k =
\delta_{\mu(i,j),k}$, and hence
  \[
  \Delta_{ij}
  =
  \sum_{k=1}^n \sum_{\ell=1}^n \delta_{\mu(j,i),k} \delta_{\mu(k,\ell),\ell}
  =
  \sum_{\ell=1}^n \delta_{\mu(\mu(j,i),\ell),\ell}.
  \]

\begin{corollary} \label{drazincorollary}
\emph{(Drazin \cite{Drazin})}
For a semigroup algebra $A$ we have
  \[
  \Delta_{ij} =   | \, \{ \, \ell \mid \mu(\mu(j,i),\ell) = \ell \, \} \, |.
  \]
\end{corollary}

To calculate a basis for the radical $R(A)$, we compute the row canonical form
$\mathrm{RCF}(\Delta)$ and extract the canonical basis for the nullspace in the
usual way: Suppose that $\Delta$ has rank $r$ and that the leading 1 of row $i$
of the RCF occurs in column $j_i$ where $1 \le j_1 < \cdots < j_r \le n$. Let
$\Lambda = \{ j_1, \dots, j_r \}$ and set $\Phi = X_n \setminus \Lambda$. For
each $k = 1, \dots, n{-}r$ set the $n{-}r$ free variables $x_j$ ($j \in \Phi$)
equal to the $k$-th unit vector in $F^{n-r}$ and solve for the leading
variables $x_j$ ($j \in \Lambda$). We obtain $n{-}r$ vectors in $F^n$ which
form a basis of $R(A)$.

\subsection{Structure constants for the semisimple quotient}

Let $\Sigma$ be the $(n{-}r) \times n$ matrix in which row $k$ contains the
coefficients of the $k$-th radical basis vector. In $\mathrm{RCF}(\Sigma)$, let
$\ell_i$ be the column containing the leading 1 of row $i$. Row $k$ of
$\mathrm{RCF}(\Sigma)$ contains the coefficients of the $k$-th reduced radical
basis vector. Set
  \[
  L = \{ \ell_1, \dots, \ell_{n-r} \},
  \qquad
  M = \{ 1, \dots, n \} \setminus L = \{ m_1, \dots, m_r \}.
  \]
The reduced radical basis vectors have the following form for some $\rho_{ij}
\in F$:
  \[
  a_{\ell_i} + \sum_{j \in M, \, j > \ell_i} \rho_{ij} a_j
  \quad
  (1 \le i \le n{-}r).
  \]
We use this reduced basis to compute the structure constants for $Q = A / R$. A
basis of $Q$ consists of the cosets $\overline{a}_m = a_m + R$ for $m \in M$.
To compute $\overline{a}_i \overline{a}_j$ we observe that $\overline{a}_i
\overline{a}_j = \overline{a_i a_j}$, but $a_i a_j$ may contain $a_\ell$ with
$\ell \in L$. These terms must be rewritten using the reduced radical basis
relations:
  \[
  \overline{a}_{\ell_i}
  =
  - \sum_{m \in M, \, m > \ell_i} \rho_{im} \overline{a}_m
  =
  - \sum_{k=1}^r \sigma_{ik} \overline{a}_{m_k},
  \quad
  \sigma_{ik}
  =
  \begin{cases}
  0 &\text{if $m_k < \ell_i$}, \\
  \rho_{im_k} &\text{if $m_k > \ell_i$}.
  \end{cases}
  \]
Because we are using the reduced basis, only $\overline{a}_m$ for $m \in M$
occur in $\overline{a}_i \overline{a}_j$. At this point we reindex the basis of
$Q$: we set $b_i = \overline{a}_{m_i}$ for $1 \le i \le r$. We have
  \allowdisplaybreaks
  \begin{align*}
  b_i b_j
  &=
  \overline{a}_{m_i} \overline{a}_{m_j}
  =
  \sum_{k=1}^n c_{m_i m_j}^k \overline{a}_k
  =
  \sum_{k = 1}^r c_{m_i m_j}^{m_k} \overline{a}_{m_k}
  +
  \sum_{h = 1}^{n-r} c_{m_i m_j}^{\ell_h} \overline{a}_{\ell_h}
  \\
  &=
  \sum_{k = 1}^r c_{m_i m_j}^{m_k} b_k
  -
  \sum_{h = 1}^{n-r} c_{m_i m_j}^{\ell_h}
  \sum_{k=1}^r \sigma_{hk} b_k
  =
  \sum_{k = 1}^r
  \Big(
  c_{m_i m_j}^{m_k}
  -
  \sum_{h = 1}^{n-r} c_{m_i m_j}^{\ell_h} \sigma_{hk}
  \Big)
  b_k.
  \end{align*}
These structure constants for $Q$ have the following form for some $d_{ij}^k
\in F$:
  \[
  b_i b_j
  =
  \sum_{k=1}^r
  d_{ij}^k b_k
  \quad
  (i, j = 1, \dots, r).
  \]

\subsection{The center of a semisimple algebra}

The next step is to compute the center $Z(Q) = \{ x \in Q \mid \text{$xy = yx$
for all $y \in Q$} \}$. We quote the following facts:

\begin{theorem}
\cite[Corollary 2.2.8, Theorem 2.4.1]{DrozdK} The center of a semisimple
algebra is semisimple. Every commutative semisimple algebra is a direct sum of
fields.
\end{theorem}

Since $Q$ is the direct sum of simple matrix algebras, and since the center of
a simple matrix algebra consists of the scalar matrices, the decomposition
  \[
  Q
  =
  Q_1 \oplus \cdots \oplus Q_c
  =
  M_{n_1}(D_1) \oplus \cdots \oplus M_{n_c}(D_c),
  \]
implies the decomposition $Z(Q) = F_1 \oplus \cdots \oplus F_c$ where $F_1,
\dots, F_c$ are extension fields of $F$. Furthermore, $Q_i = Q F_i$ for $1 \le
i \le k$, and this reduces the problem to the commutative case: if we can
decompose $Z(Q)$ into the direct sum of fields, then we can decompose $Q$ into
the direct sum of simple matrix algebras.

We can represent $Z(Q)$ as the nullspace of a matrix. Let $b_1, \dots, b_r$ be
a basis of $Q$ with structure constants $d_{ij}^k$. Then $x \in Z(Q)$ if and
only if $x b_i = b_i x$ for $1 \le i \le r$. We have
  \begin{align*}
  &
  x = \sum_{j=1}^r x_j b_j,
  \quad
  b_i x
  =
  \sum_{j=1}^r x_j b_i b_j
  =
  \sum_{j=1}^r x_j \sum_{k=1}^r d_{ij}^k b_k
  =
  \sum_{k=1}^r \sum_{j=1}^r d_{ij}^k x_j b_k,
  \\
  &
  x b_i
  =
  \sum_{k=1}^r \sum_{j=1}^r d_{ji}^k x_j b_k,
  \quad
  b_i x - x b_i
  =
  \sum_{k=1}^r \bigg( \sum_{j=1}^r (d_{ij}^k-d_{ji}^k) x_j \bigg) b_k.
  \end{align*}

\begin{corollary}
The center $Z(Q)$ is the nullspace of the $r^2 \times r$ matrix in which the
entry in row $(i{-}1)r + k$ and column $j$ is $d_{ij}^k - d_{ji}^k$ for $1 \le
i, j, k \le r$.
\end{corollary}

We compute the RCF of this matrix and the canonical basis of row vectors $z_1,
\dots, z_c$ for the nullspace. For $1 \le i, j \le c$ we use the structure
constants for $Q$ to compute a row vector $v_{ij}$ representing $z_i z_j$ as a
linear combination of $b_1, \dots, b_r$. The coefficients of $z_i z_j$ with
respect to the basis $z_1, \dots, z_c$ are the first $c$ entries in the last
column of the RCF of the following augmented matrix:
  \[
  \left[
  \begin{array}{cccc}
  z_1^t & \cdots & z_c^t & v_{ij}^t
  \end{array}
  \right].
  \]
From this we obtain the structure constants for $Z(Q)$ where $f_{ij}^k \in F$:
  \[
  z_i z_j = \sum_{k=1}^c f_{ij}^k z_k
  \quad
  (1 \le i, j \le c).
  \]

\subsection{Orthogonal idempotents in a commutative semisimple algebra}

Our next task is to decompose the commutative semisimple algebra $Z = Z(Q)$
into a direct sum of fields. We need to find a new basis $e_1, \dots, e_c$ of
orthogonal primitive idempotents: $e_i^2 = e_i$ and $e_i e_j = 0$ ($i \ne j$).
We use a recursive ideal-splitting procedure following Ivanyos and R\'onyai
\cite{IvanyosR}. Let $u$ be a (nonzero) element of a commutative semisimple
algebra $Z$. We compute a basis for the ideal $I$ generated by $u$, and
calculate the identity element of $I$. We choose a basis element $v$ of $I$
that is not a scalar multiple of the identity element. We compute the minimal
polynomial $f$ of $v$ as an element of $I$, and factor $f$ over $F$. We have
two cases:
  \begin{enumerate}
  \item[(a)] If $f$ is irreducible, then $F(v)$ is a field. If $F(v) = I$
      then we are done: the ideal $I$ is a field. If $F(v) \ne I$ then we
      choose a basis element $w$ of $I$ with $w \notin F(v)$ and compute
      the minimal polynomial of $w$ over $F(v)$. We repeat this process
      until we have either ($i$) constructed a proof that $I$ is a field or
      ($ii$) found an element of $I$ whose minimal polynomial is reducible
      over $F$.
  \item[(b)] If $f$ is reducible, then $f = gh$ where $g, h \in F[x]
      \setminus F$ are relatively prime. Hence there exist $s, t \in F[x]$
      for which $sg + th = 1$. It follows that the ideals $J$ and $K$
      generated by $g(v)$ and $h(v)$ split $I$: that is, $J$ and $K$ are
      proper ideals of $I$ such that $I = J \oplus K$ and $JK = \{0\}$.
  \end{enumerate}
This algorithm starts with $I = Z(Q)$ and recursively performs (a) and (b) to
decompose $Z$ into the direct sum of fields. It uses three subprocedures: (1)
given a generator of an ideal $I$, compute a basis of $I$; (2) given a basis of
$I$, compute the identity element of $I$; (3) given an element of $I$, compute
its minimal polynomial.

For subprocedure (1), we start with an element $u \in Z$. We use the structure
constants of $Z$ to compute the products $z_i u$ for $i = 1, \dots, c$. We put
these products into the rows of a $c \times c$ matrix and compute its RCF. The
nonzero rows of the RCF form a basis of the ideal $I$ generated by $u$.

For subprocedure (2), let $z_1, \dots, z_c$ be a basis of $Z$ and let $I$ be an
ideal with basis $y_1, \dots, y_d$. We consider an arbitrary $x \in I$ and
express $y_j$ in terms of $z_k$:
  \begin{align*}
  x y_k
  &=
  \Big( \sum_{j=1}^d x_j y_j \Big) y_k
  =
  \sum_{j=1}^d x_j y_j y_k
  =
  \sum_{j=1}^d
  x_j
  \sum_{\ell=1}^c y_{j\ell} z_\ell
  \sum_{m=1}^c y_{km} z_m
  \\
  &=
  \sum_{j=1}^d
  \sum_{\ell=1}^c
  \sum_{m=1}^c
  x_j y_{j\ell} y_{km}
  ( z_\ell z_m )
  =
  \sum_{j=1}^d
  \sum_{\ell=1}^c
  \sum_{m=1}^c
  x_j y_{j\ell} y_{km}
  \sum_{p=1}^c
  f_{\ell m}^p z_p
  \\
  &=
  \sum_{p=1}^c
  \bigg(
  \sum_{j=1}^d
  \Big(
  \sum_{\ell=1}^c
  \sum_{m=1}^c
  y_{j\ell} y_{km}
  f_{\ell m}^p
  \Big)
  x_j
  \bigg)
  z_p.
  \end{align*}
The conditions $x y_k = y_k$ for $1 \le k \le d$ give a linear system of $cd$
equations in the $d$ variables $x_1, \dots, x_d$:
  \[
  \sum_{j=1}^d
  \Big(
  \sum_{\ell=1}^c
  \sum_{m=1}^c
  y_{j\ell} y_{km}
  f_{\ell m}^p
  \Big)
  x_j
  =
  y_{kp}
  \quad
  (1 \le k \le d,
  \;
  1 \le p \le c).
  \]
The unique solution of this system is the identity element $e$ of the ideal
$I$.

For subprocedure (3), we start with an element $u \in I$, and the previously
computed identity element $e \in I$. We represent $e$ as a column vector with
respect to the basis $z_1, \dots, z_c$. Assume that for $j \ge 1$ we have
already computed the $c \times j$ matrix whose column vectors are $u^{j-1},
\dots, u, e$ and that this matrix has rank $j$; this holds when $j = 1$. We use
the structure constants for $Z$ to multiply the first column by $u$, obtaining
$u^j$; we then augment the matrix on the left. If this $c \times (j{+}1)$
matrix has rank $j{+}1$, we repeat; otherwise, we have a dependence relation
among $u^j, \dots, u, e$, and this is the (not necessarily monic) minimal
polynomial. The coefficients of the minimal polynomial are the last column of
the RCF.

\subsection{Bases for the simple ideals of the semisimple quotient}

We now have a new basis $e_1, \dots, e_c$ of orthogonal idempotents in $Z(Q)$;
these elements are the identity elements in the extension fields in the
decomposition $Z(Q) = F_1 \oplus \cdots \oplus F_c$; and these fields are the
centers of the simple ideals $Q_i = M_{n_i}(D_i)$ in the decomposition $Q = Q_1
\oplus \cdots \oplus Q_c$. We have the coefficients of $e_1, \dots, e_c$ with
respect to the basis $z_1, \dots, z_c$ of $Z(Q)$, and the coefficients of $z_1,
\dots, z_c$ with respect to the basis $b_1, \dots, b_r$ of $Q$. We obtain
elements $e_i \in Q$ (note the ambiguous notation):
  \[
  e_i
  =
  \sum_{j=1}^c e_{ij} z_j
  =
  \sum_{j=1}^c e_{ij} \sum_{k=1}^r z_{jk} b_k
  =
  \sum_{k=1}^r \Big( \sum_{j=1}^c e_{ij} z_{jk} \Big) b_k.
  \]
These elements of $Q$ are the identity matrices in the matrix algebras $Q_i =
M_{n_i}(D_i)$; they are orthogonal idempotents in $Q$, but $e_i$ is primitive
if and only if $n_i = 1$. We compute a basis of $Q_i$ by constructing a $2r
\times r$ matrix; in row $j$ of the upper (resp.~lower) $r \times r$ block we
put the coefficients of $b_j e_i$ (resp.~$e_i b_j$), with respect to $b_1,
\dots, b_r$. We compute the RCF; the nonzero rows form a basis of $Q_i$.

\subsection{Isomorphism of a simple ideal with a full matrix algebra }

Suppose that we have a basis $s_1, \dots, s_{q^2}$ and structure constants for
an algebra $S$ isomorphic to $M_q(F)$. To construct an explicit isomorphism, we
need to find a new basis $E_{ij}$ ($1 \le i, j \le q$) satisfying the matrix
unit relations $E_{ij} E_{k\ell} = \delta_{jk} E_{i\ell}$. This is easy if we
can find a basis for a minimal ($q$-dimensional) left ideal $I \subset S$: we
identify the basis elements of $I$ with the standard basis $U_1, \dots, U_q \in
F^q$, and solve the linear equations $E_{ij} U_k = \delta_{jk} U_i$ to
determine the elements $E_{ij}$. If $F$ is finite, then this can be done in
polynomial time; but if $F = \mathbb{Q}$, then the problem is more difficult,
and seems to be equivalent to hard number-theoretic problems such as integer
factorization; see R\'onyai \cite{Ronyai2}. If we are lucky, one of the basis
elements of $S$ generates a minimal left ideal; this happens in the example in
Part 2.

\subsection{Explicit matrices for the irreducible representations}

Suppose that we have found an explicit isomorphism of each simple ideal with a
full matrix algebra. We then have a new basis of $Q = Q_1 \oplus \cdots \oplus
Q_c$ consisting of matrix units:
  \[
  E_{ij}^{(k)} \in Q_k \approx M_{q_k}(F)
  \quad
  (1 \le k \le c,
  \;
  1 \le i, j \le q_k).
  \]
Let $M$ be the $r \times r$ matrix which expresses the matrix units
$E_{ij}^{(k)}$, ordered in some way, in terms of the original basis: the
$(\ell,m)$ entry of $M$ is the coefficient of $b_\ell$ in the $m$-th matrix
unit. The inverse matrix expresses the original basis in terms of the matrix
units, and has a horizontal block structure: for each $k = 1, \dots, c$ the
rows of $M^{-1}$ with indices $m$ from $q_1^2 {+} \cdots {+} q_{k-1}^2 {+} 1$
to $q_1^2 {+} \cdots {+} q_k^2$ define the projection of $Q$ onto $Q_k$. The
$\ell$-th column of the $k$-th horizontal block contains the matrix entries in
the projection of $b_\ell$ onto $M_{q_k}(F)$, and from this we obtain the
matrix for $b_\ell$ in the $k$-th irreducible representation. Composing the map
$A \to A/R = Q$ with the projection $Q \to Q_k$ gives the matrices representing
the basis elements of $A$.

\subsection{Lifting the semisimple quotient to a subalgebra}

The last step is to find a subalgebra $B \subseteq A$ which is isomorphic to
the semisimple quotient $Q$ and is a vector space complement to the radical
$R$; the existence of $B$ is guaranteed by the Wedderburn-Malcev Theorem. Let
$A$ be an associative algebra of dimension $n$ over $F$ with radical $R$ and
semisimple quotient $Q = A / R$. Let $\overline{\beta}_1, \dots,
\overline{\beta}_r$ be a basis of $Q$ where $\overline{\beta}_i = \beta_i + R$
with $\beta_i \in A$. We need to find $\gamma_1, \dots, \gamma_r \in R$ so that
$\beta_i + \gamma_i \in A$ have the same structure constants $d_{ij}^k \in F$
as $\overline{\beta}_i \in A / R$; that is,
  \[
  ( \beta_i  +  \gamma_i ) ( \beta_j  +  \gamma_j )
  =
  \sum_{k=1}^r
  d_{ij}^k ( \beta_k  + \gamma_k )
  \quad \text{where} \quad
  \overline{\beta}_i \overline{\beta}_j
  =
  \sum_{k=1}^r
  d_{ij}^k \overline{\beta}_k.
  \]
These equations can be rewritten as follows where $\delta_{ij} \in R$:
  \[
  \beta_i \beta_j
  =
  \sum_{k=1}^r
  d_{ij}^k \beta_k
  +
  \delta_{ij},
  \quad
  \beta_i \beta_j + \beta_i \gamma_j + \gamma_i \beta_j + \gamma_i \gamma_j
  =
  \sum_{k=1}^r
  d_{ij}^k \beta_k
  +
  \sum_{k=1}^r
  d_{ij}^k \gamma_k.
  \]
We combine these equations, and consider the special case in which $R^2 =
\{0\}$:
  \[
  \beta_i \gamma_j + \gamma_i \beta_j + \gamma_i \gamma_j - \sum_{k=1}^r
  d_{ij}^k \gamma_k
  =
  - \delta_{ij},
  \qquad
  \beta_i \gamma_j + \gamma_i \beta_j  - \sum_{k=1}^r d_{ij}^k \gamma_k
  =
  - \delta_{ij}.
  \]
The last equation is a linear system in the coefficients $x_{i\ell} \in F$ of
the radical terms $\gamma_i$ with respect to a basis $\zeta_1, \dots,
\zeta_{n-r}$ of $R$. We have
  \[
  \gamma_i = \sum_{\ell=1}^{n-r} x_{i\ell} \zeta_\ell,
  \qquad
  \sum_{\ell=1}^{n-r} x_{j\ell} \beta_i \zeta_\ell
  +
  \sum_{\ell=1}^{n-r} x_{i\ell} \zeta_\ell \beta_j
  -
  \sum_{k=1}^r
  d_{ij}^k
  \sum_{\ell=1}^{n-r} x_{k\ell} \zeta_\ell
  =
  - \delta_{ij}.
  \]
We expand $\beta_i \zeta_\ell$, $\zeta_\ell \beta_j$ and $\delta_{ij}$ in terms
of $\zeta_1, \dots, \zeta_{n-r}$ where $\lambda_{i\ell}^t, \rho_{\ell j}^t,
\sigma_{ij}^t \in F$:
  \[
  \beta_i \zeta_\ell
  =
  \sum_{t=1}^{n-r} \lambda_{i\ell}^t \zeta_t,
  \qquad
  \zeta_\ell \beta_j
  =
  \sum_{t=1}^{n-r} \rho_{\ell j}^t \zeta_t,
  \qquad
  \delta_{ij}
  =
  \sum_{t=1}^r
  \sigma_{ij}^t \zeta_t.
  \]
We obtain
  \[
  \sum_{\ell=1}^{n-r} x_{j\ell} \sum_{t=1}^{n-r} \lambda_{i\ell}^t \zeta_t
  +
  \sum_{\ell=1}^{n-r} x_{i\ell} \sum_{t=1}^{n-r} \rho_{\ell j}^t \zeta_t
  -
  \sum_{t=1}^{n-r}
  \sum_{k=1}^r
  d_{ij}^k
  x_{kt} \zeta_t
  =
  -
  \sum_{t=1}^{n-r}
  \sigma_{ij}^t \zeta_t.
  \]
Extracting the coefficient of $\zeta_t$ gives
  \[
  \sum_{\ell=1}^{n-r}
  \lambda_{i\ell}^t
  x_{j\ell}
  +
  \sum_{\ell=1}^{n-r}
  \rho_{\ell j}^t
  x_{i\ell}
  -
  \sum_{k=1}^r
  d_{ij}^k
  x_{kt}
  =
  - \sigma_{ij}^t
  \quad
  (1 \le i, j \le r, 1 \le t \le n{-}r).
  \]
The terms $\gamma_i$ are the solution of these $r^2 (n{-}r)$ linear equations
in the $r (n{-}r)$ variables $x_{i\ell}$. This solution is not unique: since
any two liftings of the quotient are unipotently conjugate by an element of the
form $1 + \zeta$ where $\zeta \in R$ (Theorem \ref{wedderburnmalcev}), the
number of parameters will equal the dimension of the radical.

In the general case where $R^2 \ne \{0\}$, suppose that $R^\nu \ne \{0\}$ but
$R^{\nu+1} = \{0\}$ for some $\nu \ge 1$; the special case $R^2 = \{0\}$
corresponds to $\nu = 1$. We sketch the approach developed by de Graaf et al.
\cite{deGraafIKR} which uses induction on $\mu = 1, \dots, \nu$. The inductive
step applies the computations in the special case to compute a lifting of
$A/R^\mu$ to $A/R^{\mu+1}$ by solving a linear system in the coefficients of
terms $\gamma_i \in R^\mu/R^{\mu+1}$ using a basis for a complement of
$R^{\mu+1}$ in $R^\mu$. At the last step, when $\mu = \nu$, we have obtained a
lifting of $A/R$ to a subalgebra of $A/R^{\nu+1} = A/\{0\} = A$.


\section{Semigroups of Boolean matrices}

\subsection{Binary relations on a finite set}

Perhaps the most general associative structure is the collection of binary
relations on a set under the operation of relational composition.

\begin{definition}
Let $n$ be a positive integer and set $X_n = \{ 1, \dots, n \}$. The power set
$P(X_n^2)$ of the Cartesian square is the collection of all binary relations on
$X_n$. The natural associative operation on $P(X_n^2)$ is relational
composition:
  \[
  R \circ S
  =
  \{\,
  (i,k)
  \mid
  \text{there exists $j \in X_n$ such that $(i,j) \in R$ and $(j,k) \in S$}
  \, \}.
  \]
We represent the relation $R \in P(X_n^2)$ as the $n \times n$ zero-one matrix
$(m_{ij})$ where $m_{ij} = 1$ if and only if $(i,j) \in R$. Relational
composition corresponds to matrix multiplication using Boolean arithmetic ($1 +
1 = 1$). This structure is called the \textbf{semigroup of binary relations on
$n$ elements} and is denoted $B_n$.
\end{definition}

The most familiar subsemigroup of $B_n$ is the \textbf{symmetric group} $S_n$,
consisting of all matrices in which each row and each column has exactly one 1.
The \textbf{symmetric inverse semigroup} $SI_n$, consisting of all matrices in
which each row and each column has at most one 1, corresponds to partial
bijections between subsets of $X_n$. The \textbf{full transformation semigroup}
$FT_n$, consisting of all matrices in which each column has exactly one 1,
corresponds to functions $X_n \to X_n$. The \textbf{partial transformation
semigroup} $PT_n$, consisting of all matrices in which each column has at most
one 1, corresponds to functions from subsets of $X_n$ to $X_n$. (These four
classes are the \textbf{classical finite transformation semigroups}; see
Ganyushkin and Mazorchuk \cite{GanyushkinM}.) The \textbf{semigroup of Hall
matrices} $HM_n$ consists of all matrices $(m_{ij})$ which contain a
permutation matrix in the sense that for some $\sigma \in S_n$ we have
$m_{i,\sigma(i)} = 1$ for $i = 1, \dots, n$. The \textbf{semigroup of
quasipermutations} $QP_n$ consists of all matrices in which each column and
each row has at least one 1. (Every Hall matrix is a quasipermutation, but the
converse is false for $n \ge 3$.) These semigroups can be regarded as
generalizations of the symmetric group; their orders are given in Table
\ref{semigrouporders}.

  \begin{table}
  \[
  \begin{array}{c|c|c|c|c|c|c}
  S_n & SI_n & FT_n & PT_n & HM_n & QP_n & B_n
  \\
  n!
  &
  \displaystyle{\sum_{i=1}^n \binom{n}{i}^2 i!}
  &
  n^n
  &
  (n{+}1)^n
  &
  \text{open \cite{EverettS}}
  &
  \displaystyle{\sum_{k=0}^{n-1} (-1)^k \binom{n}{k} ( 2^{n-k}{-}1 )^n}
  &
  2^{n^2}
  \end{array}
  \]
  \caption{Orders of subsemigroups of the semigroup of binary relations}
  \label{semigrouporders}
  \end{table}

For the semigroup algebra of a finite semigroup, we have two different bases:
first, the elements of the semigroup; second, the matrix units in the
Wedderburn decomposition together with the reduced basis of the radical. The
projections onto the simple ideals in the semisimple quotient provide
irreducible representations of the semigroup; see Bremner and El Bachraoui
\cite{BremnerE} for a general result regarding $B_n$.

\subsection{The partial transformation semigroup on two elements}

We explicitly compute the structure of the semigroup algebra $A = \mathbb{Q}
PT_2$ of the semigroup $\{ a_1, \dots, a_9 \}$ of all $2 \times 2$ zero-one
matrices in which each column has at most one 1:
  \[
  \left\{
  \begin{bmatrix} 0 & 0 \\ 0 & 0 \end{bmatrix},
  \begin{bmatrix} 1 & 0 \\ 0 & 0 \end{bmatrix},
  \begin{bmatrix} 0 & 1 \\ 0 & 0 \end{bmatrix},
  \begin{bmatrix} 0 & 0 \\ 1 & 0 \end{bmatrix},
  \begin{bmatrix} 0 & 0 \\ 0 & 1 \end{bmatrix},
  \begin{bmatrix} 1 & 1 \\ 0 & 0 \end{bmatrix},
  \begin{bmatrix} 1 & 0 \\ 0 & 1 \end{bmatrix},
  \begin{bmatrix} 0 & 1 \\ 1 & 0 \end{bmatrix},
  \begin{bmatrix} 0 & 0 \\ 1 & 1 \end{bmatrix}
  \right\}.
  \]
The multiplication in $PT_2$ is displayed in Table \ref{multiplicationtable}:
$a_i a_j = a_{\mu(i,j)}$ where $\mu(i,j)$ is the entry in row $i$ and column
$j$. We study $A$ since ($i$) it is a small algebra with a nonzero radical in
characteristic 0; ($ii$) there is a unique irreducible representation of
dimension $> 1$; ($iii$) the minimal polynomials of the central elements have
rational roots; ($iv$) the radical has square zero, so we can lift the quotient
in one step.

  \begin{table}
  \[
  \left[
  \begin{array}{ccccccccc}
  1 & 1 & 1 & 1 & 1 & 1 & 1 & 1 & 1 \\
  1 & 2 & 3 & 1 & 1 & 6 & 2 & 3 & 1 \\
  1 & 1 & 1 & 2 & 3 & 1 & 3 & 2 & 6 \\
  1 & 4 & 5 & 1 & 1 & 9 & 4 & 5 & 1 \\
  1 & 1 & 1 & 4 & 5 & 1 & 5 & 4 & 9 \\
  1 & 2 & 3 & 2 & 3 & 6 & 6 & 6 & 6 \\
  1 & 2 & 3 & 4 & 5 & 6 & 7 & 8 & 9 \\
  1 & 4 & 5 & 2 & 3 & 9 & 8 & 7 & 6 \\
  1 & 4 & 5 & 4 & 5 & 9 & 9 & 9 & 9
  \end{array}
  \right]
  \]
  \caption{Multiplication table for $PT_2$}
  \label{multiplicationtable}
  \end{table}

  \begin{table}
  \[
  \left[
  \begin{array}{ccccccccc}
  1 & 1 & 1 & 1 & 1 & 1 & 1 & 1 & 1 \\
  1 & 4 & 1 & 1 & 1 & 4 & 4 & 1 & 1 \\
  1 & 1 & 1 & 4 & 1 & 1 & 1 & 4 & 4 \\
  1 & 1 & 4 & 1 & 1 & 4 & 1 & 4 & 1 \\
  1 & 1 & 1 & 1 & 4 & 1 & 4 & 1 & 4 \\
  1 & 4 & 1 & 4 & 1 & 4 & 4 & 4 & 4 \\
  1 & 4 & 1 & 1 & 4 & 4 & 9 & 1 & 4 \\
  1 & 1 & 4 & 4 & 1 & 4 & 1 & 9 & 4 \\
  1 & 1 & 4 & 1 & 4 & 4 & 4 & 4 & 4
  \end{array}
  \right]
  \quad
  \left[
  \begin{array}{rrrrrrrrr}
  1 & . & . & . & . &\!\!\!\! -1 & . & . &\!\!\!\! -1 \\
  . & 1 & . & . & . &\!\!\!\!  1 & . & . &\!\!\!\!  . \\
  . & . & 1 & . & . &\!\!\!\!  1 & . & . &\!\!\!\!  . \\
  . & . & . & 1 & . &\!\!\!\!  . & . & . &\!\!\!\!  1 \\
  . & . & . & . & 1 &\!\!\!\!  . & . & . &\!\!\!\!  1 \\
  . & . & . & . & . &\!\!\!\!  . & 1 & . &\!\!\!\!  . \\
  . & . & . & . & . &\!\!\!\!  . & . & 1 &\!\!\!\!  . \\
  . & . & . & . & . &\!\!\!\!  . & . & . &\!\!\!\!  . \\
  . & . & . & . & . &\!\!\!\!  . & . & . &\!\!\!\!  .
  \end{array}
  \right]
  \]
  \caption{The radical matrix for $A = \mathbb{Q} PT_2$ and its row canonical form}
  \label{radicalmatrix}
  \end{table}

  \begin{table}
  \[
  \left[
  \begin{array}{rrrrrrrrr}
  1 &\!\!\!\!  -1 &\!\!\!\!  -1 &\!\!\!\!   . &\!\!\!\!   . &  1 &  . & . &  . \\
  1 &\!\!\!\!   . &\!\!\!\!   . &\!\!\!\!  -1 &\!\!\!\!  -1 &  . &  . & . &  1
  \end{array}
  \right]
  \quad
  \left[
  \begin{array}{rrrrrrrrr}
  1 & . & . &\!\!\!\! -1 &\!\!\!\! -1 &\!\!\!\!  . & . & . & 1 \\
  . & 1 & 1 &\!\!\!\! -1 &\!\!\!\! -1 &\!\!\!\! -1 & . & . & 1
  \end{array}
  \right]
  \]
  \caption{The canonical and reduced bases of the radical of $\mathbb{Q} PT_2$}
  \label{radicalbasis}
  \end{table}

From the multiplication table we obtain the matrix $\Delta$ which has the
radical $R$ as its nullspace (Corollary \ref{drazincorollary}), and we compute
its RCF; see Table \ref{radicalmatrix}. The matrix has rank 7, and so $R$ has
dimension 2. We set the free variables $(x_6, x_9)$ equal to $(1,0)$ and
$(0,1)$ to obtain the canonical basis of the nullspace and then the reduced
basis of the radical; see Table \ref{radicalbasis}. The reduced basis consists
of these elements of $A$:
  \begin{align*}
  \zeta_1
  &=
    \begin{bmatrix} 0 & 0 \\ 0 & 0 \end{bmatrix}
  - \begin{bmatrix} 0 & 0 \\ 1 & 0 \end{bmatrix}
  - \begin{bmatrix} 0 & 0 \\ 0 & 1 \end{bmatrix}
  + \begin{bmatrix} 0 & 0 \\ 1 & 1 \end{bmatrix},
  \\
  \zeta_2
  &=
    \begin{bmatrix} 1 & 0 \\ 0 & 0 \end{bmatrix}
  + \begin{bmatrix} 0 & 1 \\ 0 & 0 \end{bmatrix}
  - \begin{bmatrix} 0 & 0 \\ 1 & 0 \end{bmatrix}
  - \begin{bmatrix} 0 & 0 \\ 0 & 1 \end{bmatrix}
  - \begin{bmatrix} 1 & 1 \\ 0 & 0 \end{bmatrix}
  + \begin{bmatrix} 0 & 0 \\ 1 & 1 \end{bmatrix}.
  \end{align*}
We have these corresponding relations in $Q = A / R$:
  \[
  \overline{a}_1 =
  \overline{a}_4 + \overline{a}_5 - \overline{a}_9,
  \qquad
  \overline{a}_2 =
  {}
  - \overline{a}_3 + \overline{a}_4 + \overline{a}_5
  + \overline{a}_6 - \overline{a}_9.
  \]
The semisimple quotient $Q$ has dimension 7. We compute the RCF of the matrix
whose nullspace is the center, and extract the canonical basis of $Z(Q)$; see
Table \ref{centerbasis}. The center has dimension 4; its structure constants
are in Table \ref{centerstructure}. We need to find a new basis of orthogonal
idempotents.

  \begin{table}
  \[
  \left[
  \begin{array}{rrrrrrr}
  1 & . & . &\!\!\!\!  1 & . & 1 & . \\
  . & 1 & . &\!\!\!\!  . & . & 1 & 1 \\
  . & . & 1 &\!\!\!\! -1 & . & . & 1
  \end{array}
  \right],
  \qquad
  \left[
  \begin{array}{rrrrrrr}
  -1 &\!\!\!\!  . &\!\!\!\!  1 & 1 & . & . & . \\
   . &\!\!\!\!  . &\!\!\!\!  . & . & 1 & . & . \\
  -1 &\!\!\!\! -1 &\!\!\!\!  . & . & . & 1 & . \\
   . &\!\!\!\! -1 &\!\!\!\! -1 & . & . & . & 1
  \end{array}
  \right].
  \]
  \caption{RCF of center matrix, and canonical center basis}
  \label{centerbasis}
  \end{table}

  \begin{table}
  \[
  \begin{array}{l|cccc}
  \cdot & z_1 & z_2 & z_3 & z_4 \\ \hline
  z_1   & z_1 & z_1 & z_4 & z_4 \\
  z_2   & z_1 & z_2 & z_3 & z_4 \\
  z_3   & z_4 & z_3 & -z_1{+}z_2{-}z_4 & -z_4 \\
  z_4   & z_4 & z_4 & -z_4 & -z_4 \\
  \end{array}
  \]
  \caption{Structure constants for $Z(Q)$}
  \label{centerstructure}
  \end{table}

To start, $I = Z(Q)$ with identity element $z_2$. Since $z_1^2 = 1$, the
minimal polynomial of $z_1$ is $f = t^2 - t$ and so we take $g = t-1$ and $h =
t$ which gives $I = J \oplus K$ where $J = \langle z_1{-}z_2 \rangle$ and $K =
\langle z_1 \rangle$. A basis for $J$ (resp.~$K$) is $z_1{-}z_2$ and
$z_3{-}z_4$ (resp.~$z_1$ and $z_4$). In $J$ the identity element is
$-z_1{+}z_2$, and $z_3{-}z_4$ has minimal polynomial $t^2 - 1$. Hence $J$
splits into 1-dimensional ideals with bases $z_1{-}z_2{+}z_3{-}z_4$ and
$z_1{-}z_2{-}z_3{+}z_4$. In $K$ the identity element is $z_1$, and $z_4$ has
minimal polynomial $t^2 + t$. Hence $K$ splits into 1-dimensional ideals with
bases $z_4$ and $z_1{+}z_4$. Scaling these basis elements so that they satisfy
the idempotent equation $e^2 = e$, we obtain
  \[
  e_1 = \tfrac12 ( -z_1{+}z_2{-}z_3{+}z_4 ),
  \quad
  e_2 = \tfrac12 ( -z_1{+}z_2{+}z_3{-}z_4 ),
  \quad
  e_3 = -z_4,
  \quad
  e_4 = z_1{+}z_4.
  \]
These primitive idempotents in $Z(Q)$ correspond to these elements of $Q$:
  \begin{alignat*}{2}
  e_1
  &=
  b_1 - b_3 - \tfrac12 b_4 + \tfrac12 b_5 - \tfrac12 b_6 + \tfrac12 b_7,
  &\qquad
  e_2
  &=
  {} - \tfrac12 b_4 + \tfrac12 b_5 + \tfrac12 b_6 - \tfrac12 b_7,
  \\
  e_3
  &=
  b_2 + b_3 - b_7,
  &\qquad
  e_4
  &=
  {} - b_1 - b_2 + b_4 + b_7.
  \end{alignat*}
The ideals in $Q$ generated by $e_1, e_2, e_3, e_4$ have dimensions 1, 1, 1, 4
and so
  \[
  Q = A/R
  \approx
  \mathbb{Q} \oplus \mathbb{Q} \oplus \mathbb{Q} \oplus M_2(\mathbb{Q}).
  \]
The 4-dimensional ideal generated by $e_4$ has basis
  \[
  \alpha = b_1-b_7, \quad
  \beta = b_2-b_7, \quad
  \gamma = b_3-b_7, \quad
  \delta = b_4-b_7.
  \]
We need to compute an explicit isomorphism of this ideal with
$M_2(\mathbb{Q})$; that is, a new basis $E_{ij}$ which satisfies the matrix
unit relations $E_{ij} E_{k\ell} = \delta_{jk} E_{i\ell}$. The dimensions of
the left ideals generated by $\alpha, \beta, \gamma, \delta$ are 4, 2, 2, 2. In
particular, $\beta$ generates a 2-dimensional left ideal with basis $U_1 = b_1
- b_3$ and $U_2 = b_2 - b_7$. We identify $U_1, U_2$ with $(1,0), (0,1) \in
\mathbb{Q}^2$, and solve for the matrix units; we obtain
  \[
  E_{11} = - b_1 + b_3 + b_4 - b_7,
  \;
  E_{12} = - b_4 + b_7,
  \;
  E_{21} = b_3 - b_7,
  \;
  E_{22} = - b_2 - b_3 + 2 b_7.
  \]
We now have two bases for $Q$: the old basis $b_1, \dots, b_7$ and the new
basis $e_1$, $e_2$, $e_3$, $E_{11}$, $E_{12}$, $E_{21}$, $E_{22}$. Let $M$ be
the matrix whose $(i,j)$ entry is the coefficient of old basis element $i$ in
new basis element $j$. The columns of $M$ express the new basis with respect to
the old basis, and hence the columns of $M^{-1}$ express the old basis with
respect to the new basis; see Table \ref{changeofbasis}.

  \begin{table}
  \begin{align*}
  M
  &=
  \left[
  \begin{array}{rrrrrrr}
   1 &\!\!\!\! 0 &\!\!\!\! 0 &\!\!\!\! -1 &\!\!\!\! 0 &\!\!\!\! 0 &\!\!\!\!  0 \\
   0 &\!\!\!\! 0 &\!\!\!\! 1 &\!\!\!\!  0 &\!\!\!\! 0 &\!\!\!\! 0 &\!\!\!\! -1 \\
  -1 &\!\!\!\! 0 &\!\!\!\! 1 &\!\!\!\!  1 &\!\!\!\! 0 &\!\!\!\! 1 &\!\!\!\! -1 \\
  -\tfrac12 &\!\!\!\! -\tfrac12 &\!\!\!\!  0 &\!\!\!\!  1 &\!\!\!\! -1 &\!\!\!\!
   0 &\!\!\!\! 0 \\[3pt]
   \tfrac12 &\!\!\!\!  \tfrac12 &\!\!\!\!  0 &\!\!\!\!  0 &\!\!\!\!  0 &\!\!\!\!
   0 &\!\!\!\! 0 \\[3pt]
  -\tfrac12 &\!\!\!\!  \tfrac12 &\!\!\!\!  0 &\!\!\!\!  0 &\!\!\!\!  0 &\!\!\!\!
   0 &\!\!\!\! 0 \\[3pt]
   \tfrac12 &\!\!\!\! -\tfrac12 &\!\!\!\! -1 &\!\!\!\! -1 &\!\!\!\!  1 &\!\!\!\!
  -1 &\!\!\!\! 2
  \end{array}
  \right],
  \quad
  M^{-1}
  =
  \left[
  \begin{array}{rrrrrrr}
   0 &\!\!\!\!  0 & 0 &\!\!\!\!  0 & 1 &\!\!\!\! -1 & 0 \\
   0 &\!\!\!\!  0 & 0 &\!\!\!\!  0 & 1 &\!\!\!\!  1 & 0 \\
   1 &\!\!\!\!  1 & 1 &\!\!\!\!  1 & 1 &\!\!\!\!  1 & 1 \\
  -1 &\!\!\!\!  0 & 0 &\!\!\!\!  0 & 1 &\!\!\!\! -1 & 0 \\
  -1 &\!\!\!\!  0 & 0 &\!\!\!\! -1 & 0 &\!\!\!\! -1 & 0 \\
   1 &\!\!\!\! -1 & 1 &\!\!\!\!  0 & 0 &\!\!\!\!  0 & 0 \\
   1 &\!\!\!\!  0 & 1 &\!\!\!\!  1 & 1 &\!\!\!\!  1 & 1
  \end{array}
  \right].
  \end{align*}
  \caption{Change of basis matrices for $Q$}
  \label{changeofbasis}
  \end{table}

Semigroup elements $a_1$, $a_2$ are congruent modulo $R$ to linear combinations
of $a_3, \dots, a_9$; combining this with $M^{-1}$ we express $a_1$, $a_2$ in
terms of the matrix units. In this way we express all nine elements of $A$ in
terms of the matrix units, and this gives the four irreducible representations;
see Table \ref{representations}. The first two are the unit and sign
representations of the symmetric group; the third is the unit representation of
the semigroup; the fourth is the irreducible 2-dimensional representation.

  \begin{table}
  \begin{center}
  \[
  \begin{array}{crrrc}
  \text{element} &
  \mathbb{Q} & \mathbb{Q} & \mathbb{Q} &\quad M_2(\mathbb{Q})
  \\
  \begin{bmatrix} 0 &  0 \\ 0 &  0 \end{bmatrix} &
  \phantom{-}0 &
  \phantom{-}0 &
  \phantom{-}1 &\quad
  \left[  \begin{array}{rr} \phantom{-}0 & \phantom{-}0 \\ 0 & 0 \end{array}  \right]
  \\[8pt]
  \begin{bmatrix} 1 &  0 \\ 0 &  0 \end{bmatrix} &
  0 &
  0 &
  1 &\quad
  \left[  \begin{array}{rr} \phantom{-}1 & \phantom{-}0 \\ -1 & 0 \end{array}  \right]
  \\[8pt]
  \begin{bmatrix} 0 &  1 \\ 0 &  0 \end{bmatrix} &
  0 &
  0 &
  1 &\quad
  \left[  \begin{array}{rr} -1 & -1 \\ 1 & 1 \end{array}  \right]
  \\[8pt]
  \begin{bmatrix} 0 &  0 \\ 1 &  0 \end{bmatrix} &
  0 &
  0 &
  1 &\quad
  \left[  \begin{array}{rr} \phantom{-}0 & \phantom{-}0 \\ -1 & 0 \end{array}  \right]
  \\[8pt]
  \begin{bmatrix} 0 &  0 \\ 0 &  1 \end{bmatrix} &
  0 &
  0 &
  1 &\quad
  \left[  \begin{array}{rr} \phantom{-}0 & \phantom{-}0 \\ 1 & 1 \end{array}  \right]
  \\[8pt]
  \begin{bmatrix} 1 &  1 \\ 0 &  0 \end{bmatrix} &
  0 &
  0 &
  1 &\quad
  \left[  \begin{array}{rr} \phantom{-}0 & -1 \\ 0 & 1 \end{array}  \right]
  \\[8pt]
  \begin{bmatrix} 1 &  0 \\ 0 &  1 \end{bmatrix} &
  1 &
  1 &
  1 &\quad
  \left[  \begin{array}{rr} \phantom{-}1 & \phantom{-}0 \\ 0 & 1 \end{array}  \right]
  \\[8pt]
  \begin{bmatrix} 0 &  1 \\ 1 &  0 \end{bmatrix} &
   1 &
  -1 &
   1 &\quad
  \left[  \begin{array}{rr} -1 & -1 \\ 0 & 1 \end{array}  \right]
  \\[8pt]
  \begin{bmatrix} 0 &  0 \\ 1 &  1 \end{bmatrix} &
  0 &
  0 &
  1 &\quad
  \left[  \begin{array}{rr} \phantom{-}0 & \phantom{-}0 \\ 0 & 1 \end{array}  \right]
  \end{array}
  \]
  \end{center}
  \caption{Irreducible representations of $PT_2$}
  \label{representations}
  \end{table}

The last step is to find a subalgebra of $A$ isomorphic to the semisimple
quotient $Q = A / R$. We use the following ordered basis of $A$:
  \[
  \beta_1 = e_1, \;
  \beta_2 = e_2, \;
  \beta_3 = e_3, \;
  \beta_4 = E_{11}, \;
  \beta_5 = E_{12}, \;
  \beta_6 = E_{21}, \;
  \beta_7 = E_{22}, \;
  \zeta_1, \;
  \zeta_2.
  \]
We compute the quantities $\delta_{ij}$ and solve a linear system for the
coefficients $x_{i\ell}$ of the terms $\gamma_i$ for which $\beta_i + \gamma_i$
satisfy the structure constants for $Q$. We obtain the following matrix, in
which row $i$ gives the coefficients of $\gamma_i$ with respect to the
semigroup elements $a_1, \dots, a_9$; the free variables are $\alpha = x_{42}$,
$\beta = x_{52}$. The number of parameters is the dimension of the radical, as
expected:
  \[
  \left[
  \begin{array}{ccccccccc}
  0 & {-}\tfrac12 & {-}\tfrac12 & \tfrac12 & \tfrac12 & \tfrac12 & 0 & 0 & {-}\tfrac12
  \\
  {-}\beta & \tfrac12{-}\alpha & \tfrac12{-}\alpha & {-}\tfrac12{+}\alpha{+}\beta &
  {-}\tfrac12{+}\alpha{+}\beta & {-}\tfrac12{+}\alpha & 0 & 0 & \tfrac12{-}\alpha{-}\beta
  \\
  1 & 0 & 0 & {-}1 & {-}1 & 0 & 0 & 0 & 1
  \\
  0 & \alpha & \alpha & {-}\alpha & {-}\alpha & {-}\alpha & 0 & 0 & \alpha
  \\
  0 & \beta & \beta & {-}\beta & {-}\beta & {-}\beta & 0 & 0 & \beta
  \\
  \alpha & 0 & 0 & {-}\alpha & {-}\alpha & 0 & 0 & 0 & \alpha
  \\
  {-}1{+}\beta & 0 & 0 & 1{-}\beta & 1{-}\beta & 0 & 0 & 0 & {-}1{+}\beta
  \end{array}
  \right]
  \]
We add the terms $\gamma_i$ to the original coset representatives $\beta_i$ to
obtain a lifted basis of a subalgebra of $A$ isomorphic to $Q$. We include the
radical basis elements $\zeta_i$ to obtain a new basis of $A$. Choosing $\alpha
= 1$, $\beta = 0$ gives the basis for $A$ in Table \ref{newbasis}. We now have
the complete decomposition of the semigroup algebra $A = \mathbb{Q} PT_2$.

  \begin{table}
  \[
  \left[
  \begin{array}{rrrrrrrrr}
  0 & -\tfrac12 & \tfrac12 & \tfrac12 & -\tfrac12 &
  0 &\;\; \tfrac12 & -\tfrac12 & 0
  \\[3pt]
  0 & -\tfrac12 & -\tfrac12 & \tfrac12 & \tfrac12 &
  0 &\;\; \tfrac12 & \tfrac12 & -1
  \\[3pt]
  1 & 0 & 0 & 0 & 0 &\;\; 0 & 0 & 0 & 0
  \\
  0 & 1 & 0 & -1 & 0 &\;\; 0 & 0 & 0 & 0
  \\
  0 & 0 & 0 & 0 & 0 &\;\; -1 & 0 & 0 & 1
  \\
  1 & 0 & 0 & -1 & 0 &\;\; 0 & 0 & 0 & 0
  \\
  -1 & 0 & 0 & 0 & 0 &\;\; 0 & 0 & 0 & 1
  \\
  1 & 0 & 0 & -1 & -1 &\;\; 0 & 0 & 0 & 1
  \\
  0 & 1 & 1 & -1 & -1 &\;\; -1 & 0 & 0 & 1
  \end{array}
  \right]
  \]
  \caption{Basis for Wedderburn decomposition of $\mathbb{Q} PT_2$}
  \label{newbasis}
  \end{table}

\subsection{Further computations}

The Maple program used to decompose $\mathbb{Q} PT_2$ can also be used to
decompose $\mathbb{Q} PT_3$ and $\mathbb{Q} PT_4$. We obtain the following
results:
  \[
  \begin{array}{lrrrl}
  A & \dim A & \dim R & \dim Q &\quad \text{structure of $Q$ (matrix sizes)}
  \\
  \mathbb{Q} PT_2 & 9 & 2 & 7 &\quad 1, 1, 1, 2
  \\
  \mathbb{Q} PT_3 & 64 & 30 & 34 &\quad 1, 1, 1, 2, 3, 3, 3
  \\
  \mathbb{Q} PT_4 & 625 & 416 & 209 &\quad 1, 1, 1, 2, 3, 3, 4, 4, 4, 6, 6, 8
  \end{array}
  \]

\subsection{A constructive approach to the structure of algebras}

Since the original paper of Friedl and R\'onyai \cite{FriedlR}, there has been
much research on polynomial-time algorithms for explicit computation of the
structure of finite-dimensional associative algebras and Lie algebras. In
addition to the references already cited, the work of Eberly and Giesbrecht
\cite{Eberly1, Eberly2, EberlyG1, EberlyG2} deserves particular mention.

\subsection{Representation theory of finite semigroups}

There is a substantial literature on the structure theory of semigroup algebras
of finite semigroups; the classical reference is Clifford and Preston
\cite{CliffordP}. A recent monograph is Ganyushkin and Mazorchuk
\cite{GanyushkinM}; see also the paper \cite{GanyushkinMS}. For the symmetric
inverse semigroup, see Solomon \cite{Solomon}. For the full transformation
semigroup, see Putcha \cite{Putcha}.


\section*{Acknowledgements}

I thank Dr.~Delaram Kahrobaei of the City University of New York for the
invitation to speak in the Special Session on Groups, Computations, and
Applications at the 2010 Spring Eastern Sectional Meeting of the American
Mathematical Society (New Jersey Institute of Technology, Newark, NJ, May
22--23, 2010).

This research was partially supported by NSERC, the Natural Sciences and
Engineering Research Council of Canada.



\begin{thebibliography}{99}

\bibitem{BremnerE} \textsc{M. R. Bremner, M. El Bachraoui:} On the semigroup
    algebra of binary relations. \textit{Commun. Algebra} (to appear).

\bibitem{CliffordP} \textsc{A. H. Clifford, G. B. Preston:} \textit{The
    Algebraic Theory of Semigroups.} American Mathematical Society, 1961.

\bibitem{CohenIW} \textsc{A. M. Cohen, G. Ivanyos, G. B. Wales:} Finding the
    radical of an algebra of linear transformations. \textit{J. Pure Appl.
    Algebra} 117/118 (1997) 177--193.

\bibitem{deGraafIKR} \textsc{W. A. de Graaf, G. Ivanyos, A. K\"uronya, L.
    R\'onyai:} Computing Levi decompositions in Lie algebras. \textit{Appl.
    Algebra Engrg. Comm. Comput.} 8 (1997) 291--303.

\bibitem{Dickson} \textsc{L. E. Dickson:} \textit{Algebras and Their
    Arithmetics.} University of Chicago Press, 1923.

\bibitem{Drazin} \textsc{M. P. Drazin:} Maschke's theorem for semigroups.
    \textit{J. Algebra} 72 (1981) 269--278.

\bibitem{DrozdK} \textsc{Y. A. Drozd, V. V. Kirichenko:}
    \textit{Finite-Dimensional Algebras.} Springer, 1994.

\bibitem{Eberly1} \textsc{W. Eberly:} Decomposition of algebras over finite
    fields and number fields. \textit{Comput. Complexity} 1 (1991) 183--210.

\bibitem{Eberly2} \textsc{W. Eberly:} Decompositions of algebras over
    $\mathbb{R}$ and $\mathbb{C}$. \textit{Comput. Complexity} 1 (1991)
    211--234.

\bibitem{EberlyG1} \textsc{W. Eberly, M. Giesbrecht:} Efficient decomposition
    of associative algebras over finite fields. \textit{J. Symbolic Comput.} 29
    (2000) 441--458.

\bibitem{EberlyG2} \textsc{W. Eberly, M. Giesbrecht:} Efficient decomposition
    of separable algebras. \textit{J. Symbolic Comput.} 37 (2004) 35--81.

\bibitem{EverettS} \textsc{C. J. Everett, P. R. Stein:} The asymptotic number
    of $(0,1)$-matrices with zero permanent. \textit{Discrete Math.} 6 (1973)
    29--34.

\bibitem{FriedlR} \textsc{K. Friedl, L. R\'onyai:} Polynomial time solutions of
    some problems in computational algebra. \textit{Proceedings of the 17th
    Annual ACM Symposium on Theory of Computing.} Association for Computing
    Machinery, 1985, pages 153--162.

\bibitem{GanyushkinM} \textsc{O. Ganyushkin, V. Mazorchuk:} \textit{Classical
    Finite Transformation Semigroups: An Introduction.} Springer, 2009.

\bibitem{GanyushkinMS} \textsc{O. Ganyushkin, V. Mazorchuk, B. Steinberg:} On
    the irreducible representations of a finite semigroup. \textit{Proc. Amer.
    Math. Soc.} 137 (2009) 3585--3592.

\bibitem{IvanyosR} \textsc{G. Ivanyos, L. R\'onyai:} Computations in
    associative and Lie algebras. Chapter 5 of \textit{Some Tapas of Computer
    Algebra.} Springer, 1999, pages 91--120.

\bibitem{Parshall} \textsc{K. H. Parshall:} Joseph H. M. Wedderburn and the
    structure theory of algebras. \textit{Arch. Hist. Exact Sci.} 32 (1985)
    223--349.

\bibitem{Putcha} \textsc{M. Putcha:} Complex representations of finite monoids.
    \textit{Proc. London Math. Soc.} (3) 73 (1996) 623--641.

\bibitem{Ronyai1} \textsc{L. R\'onyai:} Computing the structure of finite
    algebras. \textit{J. Symbolic Comput.} 9 (1990) 355-373.

\bibitem{Ronyai2} \textsc{L. R\'onyai:} Simple algebras are difficult.
    \textit{Proceedings of the 19th Annual ACM Symposium on Theory of
    Computing.} Association for Computing Machinery, 1987, pages 398--408.

\bibitem{Solomon} \textsc{L. Solomon:} Representations of the rook monoid.
    \textit{J. Algebra} 256 (2002) 309--342.

\end{thebibliography}
\end{document}